\documentclass{llncs}
\usepackage{amsmath,amssymb,mathptmx}
\usepackage{graphicx}
\let\ge=\geqslant
\let\le=\leqslant

\begin{document}

\pagestyle{headings}     
\mainmatter              
\title{Kolmogorov complexity, Lovasz local lemma and critical exponents}
\titlerunning{Kolmogorov complexity, Lovasz local lemma and critical exponents}
\author{A.~Rumyantsev}
\authorrunning{A.~Rumyantsev}
\tocauthor{A.~Rumyantsev (Moscow State University)}
\institute{Moscow State University, Russia, Mathematics
Department, Logic and algorithms theory division}

\maketitle              

\begin{abstract}
D.~Krieger and J.~Shallit have proved that every real number
greater than~$1$ is a critical exponent of some
sequence~\cite{krieger-shallit}. We show how this result can be
derived from some general statements about sequences whose
subsequences have (almost) maximal Kolmogorov complexity. In this
way one can also construct a sequence that has no ``approximate''
fractional powers with exponent that exceeds a given value.
\end{abstract}

\section{Kolmogorov complexity of subsequences}

Let $\omega=\omega_0\omega_1\ldots$ be an infinite binary
sequence. For any finite set $A\subset \mathbb{N}$ let
$\omega(A)$ be a binary string of length $\#A$ formed by
$\omega_i$ with $i\in A$ (in the same order as in $\omega$). We
want to construct a sequence $\omega$ such that strings
$\omega(A)$ have high Kolmogorov complexity for all simple $A$.
(See~\cite{li-vitanyi} for the definition and properties of
Kolmogorov complexity. We use prefix complexity and denote it by
$K$, but plain complexity can also be used with minimal
changes.)

\begin{theorem}
        \label{subsequence1}
Let $\gamma$ be a positive real number less than~$1$. Then there
exists a sequence $\omega$ and an integer $N$ such that for any
finite set $A$ of cardinality at least $N$ the inequality
        $$
K(A,\omega(A)|t)\ge \gamma\cdot\#A
        $$
holds for some $t\in A$.
\end{theorem}

Here $K(A,\omega(A)|t)$ is conditional Kolmogorov complexity of
a pair $(A,\omega(A))$ relative to~$t$.

\textbf{Proof}. This result is a consequence of Lovasz local
lemma (see, e.g.,~\cite{randomized} for a proof):

\textbf{Lemma}. Assume that a finite sequence of events
$A_1,\ldots,A_n$ is given, for each $i$ some subset
$N(i)\subset\{1,\ldots, n\}$ of ``neighbors'' is fixed, positive
reals $\varepsilon_1,\ldots,\varepsilon_n$ are chosen in such a
way that
        $$
\Pr[A_i] \le \varepsilon_i \prod_{j\in N(i),j\ne i}
                                           (1-\varepsilon_j)
        $$
and for every $i$ the event $A_i$ is independent of the family
of all $A_j$ with $j\notin N(i), j\ne i$. Then the probability
of the event ``not $A_1$ and not $A_2$ and\ldots\ and not
$A_n$'' is at least $(1-\varepsilon_1)\cdot\ldots\cdot
(1-\varepsilon_n)$.

The standard compactness argument shows that it is enough (for
some $N$; the choice of $N$ will be explained later) to
construct an arbitrarily long finite sequence $\omega$ that
satisfies the statement of Theorem~\ref{subsequence1}. Let us
fix the desired length of this (long) sequence. For any set $A$
(whose elements do not exceed this length) and any string $Z$ of
length $|A|$ such that $K(A,Z|t)<\gamma\cdot\#A$ for every $t\in
A$ consider the event $\omega(A)=Z$; the set $A$ is callled the
\emph{support} of this event. We have to prove that the
complements of these events have non-empty intersection.

This is done by using Lovasz lemma. Let us choose some $\beta$
between $\gamma$ and $1$. Let $\varepsilon_i$ be $2^{-\beta s}$
where $s$ is the size of support of $i$th event. For each event
$\omega(A)=Z$ the neighbor events are events $\omega(A')=Z'$
such that the supports $A$ and $A'$ have nonempty intersection.
Let us check the assumptions of Lovasz lemma.

First, an event $A_i$ is independent of any family of events
whose supports do not intersect the support of $A_i$.

Second, let $\omega(A)=Z$ be an event and let $n$ be the
cardinality of~$A$. The probability of this event is $2^{-n}$.
We have to check that $2^{-n}$ does not exceed $2^{-\beta n}$
multiplied by the product of $(1-2^{-\beta m})$ factors for all
neighbor events (where $m$ is the size of the support of the
corresponding events).

This product can be split into parts according to possible
intersection points. (If there are several intersection points,
let us select and fix one of them.) Then for any $t\in A$ and
for any $m$ there is at most $2^{\gamma m}$ factors that belong
to the $t$-part and have size $m$, since there exist at most
$2^{\gamma m}$ objects that have complexity less than $\gamma m$
(relative to~$t$). Then we take a product over all $m$ and
multiply the results for all $t$ (there are $n$ of them). The
condition of Lovasz lemma (that we need to check) gets the form
        $$
2^{-n} \le 2^{-\beta n} \prod_{m>N} (1-2^{-\beta m})^{2^{\gamma m} n}
        $$
or (after we remove the common exponent $n$)
        $$
2^{\beta-1} \le \prod_{m>N} (1-2^{-\beta m})^{2^{\gamma m}}
        $$
Bernoulli inequality guarantees that this is true if
        $$
2^{\beta-1} \le 1-\sum_{m>N}  2^{\gamma m}2^{-\beta m}
        $$
Since the left hand side is less than~$1$ and the geometric
series converges, this inequality is true for a suitable~$N$.
(Let us repeat how the proof goes: we start with
$\beta\in(\gamma,1)$, then we choose $N$ using the convergence
of the series, then for any finite number of events we apply
Lovasz lemma, and then we use compactness.)

(End of proof)

The inequality established in this theorem has an useful
corollary:
        $$
K(\omega(A)|t)\ge \gamma\cdot\# A -K(A|t)-O(1),
        $$
since $K(A,\omega(A)|t)\le K(A|t)+K(\omega(A)|t)+O(1)$. For
example, if $A$ is an interval, then $K(A|t)$ is $o(\#A)$, so
this term (as well as an additive constant $O(1)$) can be
absorbed by a small change in $\gamma$ and we obtain the
following corollary (``Levin's lemma'', see~\cite{stacs2006} for
a discussion and further references): for any $\gamma<1$ there
exists a sequence $\omega$ such that all its substrings of
sufficiently large length $n$ have complexity at least $\gamma
n$.

\section{Critical exponents}

Let $X$ be a string over some alphabet, and let $Y$ be its
prefix. Then the string $Z=X\ldots XY$ is called a
\emph{fractional power} of $X$ and the ratio $|Z|/|X|$ is its
\emph{exponent}. A \emph{critical exponent} of an infinite
sequnce $\omega$ is the least upper bound of all exponents of
fractional powers that are substrings of $\omega$. D.~Krieger and
J.~Shallit~\cite{krieger-shallit} have proved the following
result:

\begin{theorem}
For any real $\alpha>1$ there exists an infinite sequence that
has critical exponent~$\alpha$.
\end{theorem}

Informally speaking, when constructing such a sequence, we need
to achieve two goals. First, we have to guarantee (for rational
numbers $r$ less than $\alpha$ but arbitrarily close
to~$\alpha$) that our sequence contains $r$-powers; second, we
have to guarantee that it does not contain $q$-powers for
$q>\alpha$. Each goal is easy to achieve when considered
separately. For the first one, we can just insert some $r$-power
for every rational $r<\alpha$. For the second goal we can use
the sequence with complex substrings: since every $q$-power has
complexity about $1/q$ of its length (the number of free bits in
it), Levin's sequence does not contain long $q$-powers if
$q>1/\gamma$.

The real problem is to combine these two goals: after we fix the
repetition pattern needed to ensure the first requirement (i.e.,
after decide which bits in a sequence should coincide) we need
to choose the values of the ``free'' bits in such a way that no
other (significant) repetitions arise. For that, let us first
prove some general statement about Kolmogorov complexity of
subsequences in the case when some bits are repeated.

\section{Complexity for sequences with repetitions}

Let $\sim$ be an equivalence relation on $\mathbb{N}$. We assume
that all equivalence classes are finite and the relation itself
is computable; moreover, we assume that for a given $x$ one can
effectively list the $x$'s equivalence class. This relation is
used as a repetition pattern: we consider only sequences
$\omega$ that follows $\sim$, i.e., only sequences $\omega$ such
that $\omega_i=\omega_j$ if $i\sim j$. For any set $A\subset
\mathbb{N}$ we consider the \emph{number of free bits in $A$},
i.e., the number of equivalence classes that have a non-empty
intersection with~$A$; it is denoted $\#_f A$ in the sequel.

There are countably many equivalence classes. Let us assign
natural numbers to them (say, in the increasing order of minimal
elements) and let $c(i)$ be the number of equivalence class that
contains~$i$. Then every sequence $\omega$ that follows the
repetition pattern $\sim$ has the form $\omega_i=\tau_{c(i)}$
for some function $c\colon\mathbb{N}\to\mathbb{N}$.

Now we assume that the equivalence relation $\sim$ (as
explained above) and a constant $\gamma<1$ are fixed.

\begin{theorem}
        \label{subsequence2}
There exists a sequence $\omega$ that follows the pattern~$\sim$
and an integer $N$ with the following properties: for every
finite set $A$ with $\#_f A\ge N$ there exists $t\in A$ such
that
        $$
K(\omega(A)|t) \ge \gamma\cdot\#_f A - K(A|t) - \log m(t)
        $$
where $m(t)$ is the ``multiplicity'' of $t$, i.e., the number of
bits in its equivalence class.
\end{theorem}

(Note that if all equivalence classes are singletons, then $\log
m(t)$ disappears, $\#_f A$ is the cardinality of $A$ and we get
an already mentioned corollary.)

\textbf{Proof}. Let $\omega_i=\tau_{c(i)}$ where $\tau$ is a
sequence that satisfies the statement of
Theorem~\ref{subsequence1} (with the same~$\gamma$). For any $A$
let $B$ be the set of all $c(i)$ for $i\in A$. Then $\# B=\#_f
A$. Theorem~\ref{subsequence1} guarantees that
$K(B,\tau(B)|u)\ge\gamma\cdot\#B$ for some $u\in B$. Since $u\in
B$, there exists some $t\in A$ such that $c(t)=u$. To specify
$t$ when $u$ is known, we need $\log m(t)$ bits, so $K(t|u)\le
\log m(t)+O(1)$. After $t$ is known, we need $K(A|t)$ additional
bits to specify $A$ and $K(\omega(A)|t)$ bits to specify
$\omega(A)$. Knowing $A$ and $\omega(A)$, we then reconstruct
$B$ and $\tau(B)$. Therefore,
        $$
\gamma\cdot\#B \le K(B,\tau(B)|u) \le \log m(t) + K(A|t) +
                    K(\omega(A)|t)+O(1),
        $$
which implies the desired inequality (with additional term
$O(1)$, which can be compensated by a small change in~$\gamma$).

\section{Construction}

Assume that $1 < \alpha <\beta$. First, let us show that
Theorem~\ref{subsequence2} implies the existence of a binary
sequence $\omega$ that contains fractional powers of all
rational exponents less than $\alpha$, but does not contain long
fractional powers of exponents greater than $\beta$.

To construct such a sequence, let $r_1,r_2,\ldots$ be all
rational numbers between $1$ and $\alpha$. For each
$r_i=p_i/q_i$ we ``implant'' a fractional power of exponent
$r_i$ in the sequence: we select some interval of length $p_i$
and decide that this interval should be a fractional power of
some string of length $q_i$ (and exponent $r_i$). This means
that we declare two indices in this interval equivalent if they
differ by a multiple of $q_i$. (The intervals for different $i$
are disjoint.) We call these intervals \emph{active intervals}.
We assume that distance between two active intervals is much
bigger than the lengths of these two intervals (see below why
this is useful).

\begin{figure}[h]
\begin{center}
\includegraphics{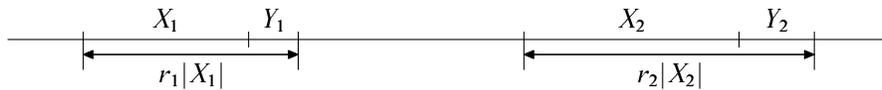}
\end{center}
\caption{Two fractional powers of exponent $r_1$ and $r_2$ are implanted;
$Y_i$ is a prefix of $X_i$ (in this example the exponents are less than $2$,
so only one full period is shown).}
\end{figure}

Evidently, any sequence that follows this repetition pattern has
critical exponent at least $\alpha$.

Let us choose some $\gamma$ between $\alpha/\beta$ and $1$ and
apply Theorem~\ref{subsequence2} with this $\gamma$ to the
pattern explained above. We get a bit sequence; let us prove
that it does not contain \emph{long} fractional powers of
exponent greater than~$\beta$. Indeed, it is easy to see that
density of free bits in this pattern is at least $1/alpha$,
i.e., for any interval $A$ of length $l$ the number of free bits
in it, $\alpha_f A$, is at least $l/\alpha$. Indeed, if $A$
intersects with two or more active intervals, then all bits
between them are free, and the distance between the intervals is
large compared to interval sizes. Then we may assume that $A$
intersects with only one active interval. All subintervals of
the active interval have the same repetitions period, and the
density of free bits is minimal when $A$ is maximal, i.e.,
coincides with the entire active interval. The bits outside the
active interval are free (no equivalences), so they can only
increase the fraction of free bits.

On the other hand, a fractional power of exponent $\beta$ and
length $l$ has complexity $l/\beta + O(\log l)$ (we specify the
length of the string and $l/\beta$ bits that form the period).
For long enough strings we then get a contradiction with the
statement of Theorem~\ref{subsequence2} since
$\alpha/\beta<\gamma$.

To get rid of short fractional powers of exponent greater
than~$\beta$ we can add additional layer of symbols that
prevents them. In other terms, consider a sequence in a finite
alphabet that follows (almost) the same repetition pattern but
has no other repetitions (not prescribed by the pattern) on
short distances. It is easy to construct such a sequence; for
example, we may assume that $q_i$ is a multiple of $i!$ and then
consider a periodic sequence with any large period $M$; it will
destroy all periods that are not multiple of $M$, i.e., all
short periods and only finitely many of $q_i$ (the latter does
not change the critical exponent). The Cartesian product of
these two sequences ($i$th letter is a pair formed by $i$th
letters of both sequences) has critical exponent between
$\alpha$ and $\beta$.

In fact, we even get a stronger result:

\begin{theorem}
        \label{sequence3}
For any $\alpha$ and $\beta$ such that $1<\alpha<\beta$ there
exist a sequence $\omega$ that has fractional powers of exponent
$r$ for all $r<\alpha$ but does not have approximate fractional
powers of exponent $\beta$ or more: there exists some
$\varepsilon>0$ such that any substring of length $n$ is
$\varepsilon n$-far from any fractional power in terms of
Hamming distance \textup(we need to change at least $\varepsilon
n$ symbols of the sequence to get a fractional power of
length~$n$\textup).
\end{theorem}

Indeed, a change of $\epsilon$-fraction bits in a sequence of
length $n$ increases its complexity at most by
$H(\varepsilon)n+O(\log n)$ where
        $$
H(\varepsilon)=-\varepsilon\log\varepsilon -
(1-\varepsilon)\log(1-\varepsilon).
        $$
Therefore, we need to change a constant fraction of bits to
compensate for the difference in complexities (between the lower
bound guaranteed by Theorem~\ref{subsequence2} and the upper
bound due to approximate periodicity). (End of proof.)

\section{Critical exponent: exact bound}

The same construction (with some refinement) can be used to
get a sequence with given critical exponent.

\begin{theorem}
\textup(Krieger -- Shallit\textup) For any real number
$\alpha>1$ there exists a sequence that has critical
exponent~$\alpha$.
\end{theorem}

(This proof follows the suggestions of D.~Krieger who informed
the author about the problem and suggested to apply
Theorem~\ref{subsequence1} to it. See~\cite{krieger-shallit} for
the original proof. Author thanks D.~Krieger for the
explanations and both authors of~\cite{krieger-shallit} for the
permission to cite their paper.)

Again, let us consider repetition pattern that guarantees all
exponents less than $\alpha$ and apply
Theorem~\ref{subsequence2} with some $\gamma$ close to $1$. This
(as we have seen) prevents powers with exponents greater that
$\alpha/\gamma$; the problem is how to get rid of intermediate
exponents.

To do this, we should distinguish between two possibilities:
(a)~an unwanted power is an extension of the prescribed one (has
the same period that unexpectedly has more repetitions) and
(b)~an unwanted power is not an extension. The first type of
unwanted powers can be prevented by adding brackets around each
active interval (in a special layer: we take a Cartesian product
of the sequence and this layer).

It remains to explain why unwanted repetitions of the second
type do not exist (for $\gamma$ close enough to~$1$). Consider
any fractional power with exponent greater than~$\alpha$. There
are two possibilities:

(1)~It intersect at least two active intervals. Then it contains
all free bits between these intervals, and (since we assume that
the distances are large compared to the length of intervals) the
density of free bits is close to $1$, so exponent greater than
$\alpha$ is impossible.

(2)~It intersects only one active interval. The same argument
(about density of free bits) shows that if the endpoints of this
fractional power deviate significantly from the endpoints of the
active interval, then the density of free bits is significantly
greater than $1/\alpha$ and we again get a contradiction.
Therefore, taking $\gamma$ close to $1$ we may guarantee that
the distance between endpoints of fractional power and active
interval is a small fraction of the length of the active
interval. Then we get \emph{two different periods} in the
intersection of fractional power and active interval. One
(``old'') is inherited from the repetition pattern; the second
one (``new'') is due to the fact that we consider a fractional
power. (The periods are different, otherwise we are in the
case~(1).) The period lengths are close to each other. Indeed,
if the new period is significantly longer, then the exponent is
less than $\alpha$; if the new period is significantly shorter,
then the complexity bound decreases and we again get a
contradiction.

Now note that two periods $t_1$ and $t_2$ in a string guarantee
the period $t_1-t_2$ near the endpoints of this string (at the
distance equal to the difference between string length and
minimal of these periods). Therefore we get a period that is a
small fraction of the string length at an interval whose length
is a non-negligible fraction of the string length. This again
significantly decreases the complexity of the string, and this
contradicts the lower bound of the complexity. (End of the
proof.)

Remark. This proof uses some parameters that have to be chosen
properly. For a given $\alpha$ we choose $\gamma$ that is close
enough to $1$ and makes the arguments about ``sufficiently small''
and ``significantly different'' things in the last paragraph valid
for long strings. Then we choose the repetition patterns where
length of active intervals are multiples of factorials and the
distances between them grow much faster than the lengths of active
intervals. Then we apply Theorem~\ref{subsequence2} for this
pattern. Finally, we look at the length $N$ provided by this
theorem and prevent all shorter periods by an additional layer.
Another layes is used for brackets. These layers destroy only
finitely many of prescribed patterns and unwanted short periods.

\end{document}